\documentclass[journal,twoside,web]{ieeecolor}

\usepackage{generic}
\usepackage{amssymb,amsmath}  

\usepackage{amsthm}

\newtheorem{assumption}{Assumption}
\newtheorem{probdef}{Problem Definition}

\usepackage{graphics} 
\usepackage{epsfig} 
\usepackage{mathptmx} 
\usepackage{times} 
\usepackage{color}
\usepackage{algorithm,algorithmic,url}
\usepackage{xfrac}
\usepackage[hidelinks]{hyperref}
\newcommand{\algrule}[1][.2pt]{\par\vskip.5\baselineskip\hrule height #1\par\vskip.5\baselineskip}
\usepackage{siunitx}

\usepackage{mdframed} 
\usepackage{xparse} 
\newcommand{\red}[1]{{\color{black}#1}}

\definecolor{darkred}{RGB}{139,0,0}
\newcommand{\Transp}{\mathsf{T}}

\newcommand{\nc}{\newcommand}
\nc{\uu}{\mathbf{u}_t}
\nc{\ww}{\mathbf{w}_t}
\nc{\xx}{\mathbf{x}_t}
\nc{\zz}{\mathbf{z}_t}
\nc{\bmat}[1]{\begin{bmatrix}#1\end{bmatrix}}
\newcommand{\pmat}[1]{\begin{bmatrix}#1 \end{bmatrix}}
\newcommand{\dx}[1]{\,\mathrm{d}#1}
\newcommand{\R}{\mathbb{R}}

\renewcommand{\baselinestretch}{0.96}

\title{\LARGE \bf
Data to Controller for Nonlinear Systems: An Approximate Solution 
}

\author{Johannes N. Hendriks, James R. Z. Holdsworth, Adrian G. Wills, Thomas B. Sch\"{o}n and Brett Ninness 
\thanks{This research was partially supported by the project \emph{NewLEADS - New Directions in Learning Dynamical Systems} (621-2016-06079), funded by the Swedish Research Council and by the \emph{Kjell och M{\"a}rta Beijer Foundation}.}
\thanks{Johannes Hendriks, James Holdsworth, Adrian Wills and Brett Ninness are with the School of Engineering, University of Newcastle, Australia.
        {\tt\small \{johannes.hendriks,james.holdsworth,\newline adrian.wills,
        brett.ninness\}@newcastle.edu.au}}%
\thanks{Thomas B. Sch\"{o}n is with the Department of Information Technology, Uppsala University, Sweden. 
        {\tt\small thomas.schon@it.uu.se}}%
}

\begin{document}
\newcommand{\coverTitle}{Data to Controller for Nonlinear Systems: An Approximate Solution}
\newcommand{\coverAuthors}{Johannes N. Hendriks, James R. Z. Holdsworth, Adrian G. Wills, Thomas B. Sch\"{o}n and Brett Ninness.}
\newcommand{\coverStatus}{Accepted for publication.}

\begin{titlepage}
    \begin{center}
        {\large \em Technical report}
        
        \vspace*{2.5cm}
        %
        {\Huge \bfseries \coverTitle  \\[0.4cm]}
        
        %
        {\Large \coverAuthors \\[2cm]}
        
        \renewcommand\labelitemi{\color{red}\large$\bullet$}
        \begin{itemize}
            \item {\Large \textbf{Please cite this version:}} \\[0.4cm]
            \large
            \coverAuthors. \coverTitle. \textit{IEEE Control Systems Letters, June 2021. doi: 10.1109/LCSYS.2021.3090349}  
        \end{itemize}
        
        \vfill

        
        \vfill
        \vspace{50mm}
    \end{center}
\end{titlepage}
\newpage
\thispagestyle{empty}
\newpage

\maketitle
\thispagestyle{empty}
\pagestyle{empty}

\begin{abstract}
This paper considers the problem of determining an optimal control
action based on observed data. We formulate the problem assuming that the system can be modelled by a nonlinear state-space model, but where the model parameters, state and future disturbances are not known and are treated as random variables. Central to our formulation is that the joint distribution of these unknown objects is conditioned on the observed data. Crucially, as new measurements become available, this joint distribution continues to evolve so that control decisions are made accounting for uncertainty as evidenced in the data. The resulting problem is intractable which we obviate by providing approximations that result in finite dimensional deterministic optimisation problems. The proposed approach is demonstrated in simulation on a nonlinear system.
\end{abstract}


\section{Introduction} 
\label{sec:introduction}

\IEEEPARstart{T}{he} control community has long been inspired by the concept of ``data to controller" -- using observed data to design a controller that not only satisfies certain performance objectives, but does so with a quantifiable level of certainty. This vision can be considered within the dual control problem setting, which dates back to the work in \cite{Feldbaum:1960}.

The dual control \red{and data-driven control problems} have recently received significant attention with a surge of new and fascinating results. \red{Data-driven focuses on control algorithms that utilise available offline and online measured data} and the dual control formulations point out the inherent trade-off between \textit{exploration}---learning a probabilistic model of the system---and \textit{exploitation}---controlling the system. 
For recent overviews of the area we \red{refer the reader} to \cite{Recht:2019,matni2019self,Mesbah:2018, hou2013model}\red{, which reveal that the literature is vast with a significant body of recent work concentrating on the linear dual control problem \cite{FerizbegovicUHS:2020,venkatasubramanian2020robust,iannelli2020multiobjective,Mesbah:2018,UmenbergerFSH:2019}, and on nonlinear data-driven control where assumptions on the uncertainty, such as Gaussian posteriors or measurement likelihoods are made \cite{lian2021koopman, yingzhao2020gaussian,carron2019data,galioto2020bayesian,khalil2015estimation}). Hence, we focus on those papers --- \cite{Mesbah:2018} and \cite{wabersich2020bayesian} ---  with a problem definition most similar to our own and highlight important differences.}


The overview in \cite{Mesbah:2018} considers \red{the nonlinear case with the crucial difference that \cite{Mesbah:2018}} deals with parameter uncertainty by assuming that the parameters are time-varying and obey the Markov property, therefore affording the use of Bayesian filtering techniques. \red{The work in \cite{wabersich2020bayesian} focuses on episodic Bayesian model predictive control (MPC), with a single parameter sample used in each episode.} In \cite{wabersich2020bayesian} it is assumed that the state can be directly measured, and sub-Gaussian distributions are assumed.

This paper addresses a more general case where: 1. The nonlinear model structure is more general than \cite{Mesbah:2018} in that we allow for a measurement model with uncertain parameters (thus also relaxing the direct state measurement assumption from \cite{wabersich2020bayesian}); 2. The requirement from \cite{Mesbah:2018} of time-varying parameters that obey the Markov property is removed by assuming a joint distribution of the parameters and the state trajectory; \red{3. Unlike \cite{wabersich2020bayesian}, we consider continuous learning rather than episodic and utilise multiple samples from the distributions.}  In addition to considering a more general setting, we also offer a practical approach for solving this problem.

The \textit{main contributions} are: 1. The formulation of a general control problem where the emphasis is moving from system data to control action; 2. Use of Hamiltonian Monte Carlo (\red{HMC}) to efficiently sample from the joint distribution of the state, parameter and future disturbances, but conditioned on all available data, resulting in a finite sum approximation to the intractable expectation integrals. Related to this we have the following more technical contributions: 3. Relaxation of so-called chance constraints, enabling direct use of smooth optimisation. This relaxation is tightened as the proposed optimisation algorithm iterates towards a solution; 4. Tailored interior-point method for solving the associated sequence of optimisation problems.

\section{Problem Definition} 
\label{sec:problem_definition}

This section details the stochastic MPC (SMPC) problem. To this end, we
assume that the system of interest has input
$u_t \in \mathbb{R}^{N_u}$ and output $y_t \in \mathbb{R}^{N_y}$ where
the integer $t$ is used to indicate a discrete time index. The system
input and output are assumed to be adequately modelled by the following
nonlinear state-space model
\begin{subequations}
	\label{eq:ssm}
	\begin{align}
    	x_{t+1} &= f(x_t, u_t, \theta, w_t), \\
	    y_t &= h(x_t, u_t, \theta, e_t),
	\end{align}
\end{subequations}
where $x \in \mathbb{R}^{N_x}$ is the state, $\theta$ are 
the parameters, and $w_t$ and $e_t$ are the process 
and measurement noise, respectively.

It is assumed that at time $t$ we have available the measured inputs
$u_{1:t} = \{u_1, \dots, u_t\}$ and measured outputs $y_{1:t}$. The 
aim is to determine a control action at time $t+1$ that minimises a user-defined 
cost. To this end it is important to discuss the presence of uncertainty 
in this problem. Therefore, in what follows we assume that
\begin{itemize}
\item \red{The structure of $f$ and $h$ is known} and the model parameters $\theta$ are treated as a \red{continuous} random variable to
  reflect uncertainty in these models;
\item Similarly, the state $x_t$ will also be treated as a random
  variable to reflect uncertainty in the model state;
\item The noise terms $w_t$ and $e_t$ are also treated as \red{continuous} random variables \red{with known, possibly non-Gaussian, distributions and} reflect state transition and measurement uncertainty. \red{These distributions can have unknown parameters to be estimated.}
\end{itemize}
Further, we assume that prior knowledge of $\theta, x_0$ can be
captured in a prior distribution $p(\theta,x_0)$: \red{for details on the choice and impact of $p(\theta,x_0)$ see \cite{Dahlin2014}.}
These random variables are conditioned on the data available \red{using an appropriate Bayesian method, such as HMC described later, to give samples from $p(x_t, \theta, \ww |y_{1:t}, u_{1:t})$.}


With this as background, we aim to determine an optimal future
input sequence that minimises some user-defined cost of a finite
future horizon length $N$, and subject to certain constraints. Towards
this, and to reduce notational overhead, we introduce the
following notation
\begin{align}
  \label{eq:vardef}
  \uu \triangleq u_{t+1:t+N+1}, \,\,\,\,
  \ww \triangleq w_{t+1:t+N+1}, \,\,\,\, 
  \xx \triangleq x_{t+1:t+N+1}.
\end{align}

\red{Importantly, once we have $\uu$ and samples of $x_t$, $\theta$ and $\ww$, then 
$\xx$ can be computed deterministically by application of~\eqref{eq:ssm}. 
That is, given samples $x_t^i, \theta^i, \ww^i$ from $p(x_t, \theta, \ww |y_{1:t}, u_{1:t})$ as well as $u_t$ and $\uu$, then samples of the future states are obtained by simulation}
\begin{subequations}
  \label{eq:state_evolve}
\begin{align}
  x^i_{t+1} &= f(x^i_t,u_t, \theta^i, w_t^i),\\
	\notag
          &\vdots\\
  x^i_{t+N+1} &= f(x_{t+N}^i,u_{t+N}, \theta^i, w^i_{t+N}).
\end{align}
\end{subequations}

\begin{probdef}
The problem of interest is
\begin{small}
\begin{align}
    \uu^* = \arg &\min_{\uu} \int V_t(x_t, u_t, \theta, \uu,\ww)\, p(x_t, \theta, \ww |y_{1:t}, u_{1:t})\, \mathrm{d}x_t \,\mathrm{d}\theta \,\mathrm{d}\ww, \notag\\
    &\text{s.t.} \quad c_u(\uu) \succeq 0,\label{eq:optimisation_problem}\\
&\qquad \mathbb{P}\left(c^j_x(x_t,u_t,\theta,\uu,\ww) \geq 0 \right) \geq 1-\epsilon, \quad j=1,\ldots,n_{cx}. \notag
\end{align}
\end{small}
\end{probdef}
Here, it is assumed that the function $c_u(\uu) \in \R^{n_{cu}}$ and $c_u(\uu) \succeq 0$ should be interpreted element-wise, and the function 
$c_x(x_t,u_t,\theta,\uu,\ww) \in \R^{n_{cx}}$ with $c^j_x(\cdot)$ indicating the $j^{\text{th}}$ element. 

The cost function $V_t(\cdot)$ is typically defined by the user to reflect a
desired outcome or performance goal. This can include both control
performance (exploitation) and possibly identification (exploration)
related goals also. \red{We make the following assumption on the cost:}
\begin{assumption}
The cost $V_t(x_t, u_t, \theta, \uu,\ww)$ is assumed to be a twice continuously differentiable
function of $\uu$ and is bounded below on the feasible domain, as
determined by the constraints. Note that this implies restrictions on
the state-space model \eqref{eq:ssm}.
\end{assumption}

The constraints $c_u(\uu) \succeq 0$ allow for modelling input
restrictions. The so-called \emph{chance constraint} (see e.g. \cite{peng2019chance,li2008chance} for further details)
$\mathbb{P}\left(c^j_x(x_t,u_t,\theta,\uu,\ww) \geq 0 \right) \geq
1-\epsilon$ is present in order to model constraints that involve random variables. The right-hand-side term $1-\epsilon$ provides the user with a mechanism to tradeoff constraint satisfaction with feasibility. In what follows, $\epsilon$ will be treated as a slack variable that will be minimised subject to a non-negativity condition. 
\red{We make the following assumption on the constraints:}
\begin{assumption}
It is assumed that the set
$\{\uu : c_u(\uu) \succeq 0 \}$ is non-empty and that $c_u(\uu)$ and $c^j_x(x_t,u_t,\theta,\uu,\ww)$ are twice continuously differentiable functions of $\uu$. 
\end{assumption}

The problem is aimed at minimising the expected value of the
cost relative to the joint conditional distribution over
$x_t, \theta$ and $\ww$ given the past measurements $u_{1:t}$ and
$y_{1:t}$. Importantly, the resulting expected value is a function of
the future control actions $\uu$ and the past data only. It is subtle, but important, to notice that this step involves full \emph{Bayesian nonlinear system identification}, which is perhaps most clear by noticing that a marginal of the required joint distribution $p(x_t, \theta, \ww |y_{1:t}, u_{1:t})$ is $p(\theta |y_{1:t}, u_{1:t})$.
That is, the posterior of the model parameters given the data. This encapsulates the full probabilistic information of the chosen model parameters available from the data. While this problem has been formulated for a long time---see e.g. \cite{Peterka1981} for an early formulation---it is only very recently (see e.g. \cite{AndrieuDoucetHolenstein2010,hendriks2020practical,LindstenJS:2014}) that we can actually solve this problem in a satisfactory way.
\section{Method} 
\label{sec:method}

There are at least two major challenges in solving problem
\eqref{eq:optimisation_problem}. Firstly, the expectation integral in
\eqref{eq:optimisation_problem} is generally intractable. Secondly,
the chance constraint is problematic since evaluation of the required
probability measure is prohibitive in all but relatively simple
cases. We are therefore faced with the following dichotomy, either
restrict the problem so that the expectations and chance constraints
are tractable, or employ an approximation that handles the
problem more generally.

We opt for the second approach and employ an approximation. 
Specifically, we will make use of a Monte Carlo 
approximation that will be applicable to both the expectation and
chance constraint. Towards this end, Section~\ref{sec:cost_approx}
will introduce the cost approximation,
Section~\ref{sec:constraint_approx} will discuss the chance constraint
approximation, and Section~\ref{sub:solving_the_optimisation_problem}
will combine these into a tractable optimisation problem. This latter
section will also provide details of an interior-point approach for
solving this optimisation problem. Finally,
Section~\ref{sub:on_line_sampling_and_control} presents an overall
SMPC algorithm, which will be demonstrated by simulation in
Section~\ref{sec:sims}.

\subsection{Cost Approximation} 
\label{sec:cost_approx}
We will now explain how the cost in~\eqref{eq:optimisation_problem} is 
replaced by a finite sum via a Monte Carlo approximation. 
Assume for the moment that we can draw $M>0$ samples from the
distribution $p(x_t, \theta, \ww|y_{1:t}, u_{1:t})$ such that
\begin{equation}
  \begin{split}
    x_t^i, \theta^i,\ww^i \sim p(x_t, \theta, \ww|y_{1:t}, u_{1:t}),
  \end{split}
\end{equation}
where the superscript $i$ will be used to indicate that this was the
$i^{\text{th}}$ such sample out of $M$. Then via the Law of Large Numbers, we
arrive at the following Monte Carlo approximation
\begin{equation}
  \begin{split}
    &\int V_t(x_t, u_t, \theta, \uu,\ww)p(x_t, \theta, \ww|y_{1:t},
    u_{1:t})\, \mathrm{d}x_t \,\mathrm{d}\theta \\ &\qquad \approx
    \frac{1}{M}\sum_{i=1}^M V_t(x_t^i, u_t, \theta^i, \uu,\ww^i).
  \end{split}
\end{equation}
The primary benefit is that this is now a finite sum and a
deterministic function of $\uu$ since the samples are fixed. \red{Under mild regularity conditions, the approximation converges almost surely to the desired expectation as $M \to \infty$~\cite{tierney1994markov}. In practice we are forced to choose a finite $M$ and deal with any consequences, which will be highly problem dependent.}   

Returning to the assumption that we are able to draw samples from the
desired conditional distribution, this is by no means trivial and
significant research attention has been directed towards solving it
(see e.g. \cite{ninness2000strong,hendriks2020practical,AndrieuDoucetHolenstein2010}). One of the most successful approaches to
producing these samples is to construct a Markov chain whose
stationary distribution coincides with the desired target distribution
$p(x_t, \theta, \ww|y_{1:t}, u_{1:t})$. Together, the use of a Markov
chain with the Monte Carlo approximation leads to the MCMC
approaches.

In this paper, \red{we suggest the use of HMC, which uses a Hamiltonian system to construct such a Markov chain. Alternates include Particle  Markov  chain  Monte Carlo \cite{AndrieuDoucetHolenstein2010} and Particle Gibbs with ancestor sampling \cite{LindstenJS:2014}}.  While the intricate details of HMC are beyond the
scope of this paper, in essence, the HMC approach is a version of the
Metropolis-Hastings algorithm~\cite{Betancourt2017,neal2010}, where the essential
steps in each iteration are described in Algorithm~\ref{alg:hmc}. \red{A discussion of the methods for determining the integration time $L$ and the mass matrix $\Sigma_k$ as well as the optimal target acceptance rate are beyond the scope of this paper and the interested reader is referred to \cite{Betancourt2017}. Further, the required derivatives are typically computed using automatic differentiation}.

\begin{algorithm}[!ht]
  \caption{\textsf{Iteration of Hamiltonian Monte Carlo}}
  \footnotesize
  \begin{algorithmic}[1]
    \STATE Let $\eta$ be the position variable of a Hamiltonian system with corresponding momentum variables $\rho$ and the current sample be given by $\eta_0 = \{x^i_{1:t},\theta^i\}$.
  Additionally, let the new target distribution for this augmented
  system be the exponential of the negative total energy,
  $\tilde{\pi}(\eta,\rho) = \exp(-H(\eta,\rho))$, where the
  Hamiltonian function is 
  \begin{align*}
    H(\eta,\rho) &= -\log p( \eta | y_{1:t}, u_{1:t}) -\log \mathcal{N}(\rho | 0,\Sigma_K),
  \end{align*}
  and $\log\mathcal{N}(\rho|0,\Sigma_k)$ is called a Euclidean-Gaussian kinetic energy with mass matrix $\Sigma_K$.
  \STATE Sample the momentum variables via $\rho_0\sim\mathcal{N}(0,\Sigma_K)$.
  \STATE Simulate the system for $L>0$ time by solving (typically using symplectic integration),
  \begin{equation}
  \pmat{\eta_L \\ \rho_L} = \int_0^L \pmat{-\frac{\partial H(\eta(\tau),\rho(\tau))}{\partial \rho} \\ \frac{\partial H(\eta(\tau),\rho(\tau))}{\partial \eta}} \dx{\tau};
  \end{equation}
\STATE Compute the Metropolis acceptance probability according to 
\begin{equation}
\begin{split}
  \alpha(\eta_L,\rho_L,\eta_0,\rho_0) &=
  \min\left\{1,\frac{\tilde{\pi}(\eta_L,-\rho_L)}{\tilde{\pi}(\eta_0,\rho_0)}\right\}\\
  &= \min\left\{1,\exp(-H(\eta_L,-\rho_L)+H(\eta_0,\rho_0))\right\}.\notag
\end{split}
\end{equation}
\STATE With probability $\alpha$, set $\{x^{i+1}_{1:t},\theta^{i+1} \}=
  \eta_L$, otherwise $\{x^{i+1}_{1:t},\theta^{i+1}
  \}= \eta_0$.
  \end{algorithmic}
  \label{alg:hmc}
\end{algorithm}

The benefit of using HMC is that it is very efficient and the
algorithm progresses with high acceptance ratio\footnote{This arises since Hamiltonian
  systems are energy preserving and hence, theoretically,
  $H(\eta_L,\rho_L)=H(\eta_0,\rho_0)$ and the acceptance ratio will be
  one.} and with low 
correlation in the chain \cite{neal2010}. For the simulations provided in this paper, we used the
\texttt{stan} software package~\cite{Stan2017} with full details provided
in the available github repository\footnote{
  \url{https://github.com/jnh277/data_to_mpc}}. Further details of the
HMC method for the current context are provided in
\cite{hendriks2020practical}.




\subsection{Chance Constraint Approximation} 
\label{sec:constraint_approx}
The chance constraint can be converted into an expectation using the
standard identity that
{\small
\begin{align}
  &\mathbb{P}\left(c^j_x(x_t, u_t,\theta,\uu,\ww) \geq 0 \right)
    \\ &\quad= \int
  \mathbb{I}(c^j_x(x_t^i, u_t,\theta^i,\uu,\ww^i)\geq 0) \, p(x_t,
  \theta, \ww |y_{1:t}, u_{1:t})\, \mathrm{d}x_t \,\mathrm{d}\theta
  \,\mathrm{d}\ww. \notag
\end{align}
}
In the above, $\mathbb{I}(\cdot)$ is the indicator function that is
one if the argument is true and zero otherwise. We can reuse the
HMC method described above to also approximate
this expectation using a finite sum via
\small
\begin{equation}
  \begin{split}
    \mathbb{P}\left(c^j_x(x_t, u_t,\theta,\uu,\ww) \hspace{-0.5mm}\geq \hspace{-0.5mm}0 \right) \hspace{-0.75mm}\approx\hspace{-0.75mm} \frac{1}{M}\sum_{i=1}^M \mathbb{I}(c^j_x(x_t^i, u_t,\theta^i,\uu,\ww^i)\geq 0),
  \end{split}
\end{equation}
\normalsize
where the samples $x_t^i,\theta^i,\ww^i$ are reused from the
cost approximation. 

Unfortunately, the indicator function is discontinuous and this would
result in a mixed-integer nonlinear programming problem, which is 
typically more challenging to solve. To avoid this, the 
indicator function is instead replaced by the logistic function
\begin{equation}
  \mathbb{I}(z\geq0) \approx \sigma(z,\gamma) = \frac{1}{1 - \exp(-\sfrac{z}{\gamma})}
\end{equation}
which coincides with the indicator function in the limit as
$\gamma \to 0$. The probability of satisfying the constraint can
therefore be approximated as
\small
\begin{equation}
  \begin{split}
    \mathbb{P}\left(c^j_x(x_t,\theta,u_{t},\uu,\ww) \hspace{-0.5mm}\geq\hspace{-0.5mm} 0 \right) &\hspace{-0.75mm}\approx\hspace{-0.75mm} \frac{1}{M}\sum_{i=1}^M \sigma(c^j_x(x_t^i,\theta^i,u_{t},\uu,\ww^i,),\gamma).
  \end{split}
\end{equation}
\normalsize

\subsection{Solving the Optimisation Problem} 
\label{sub:solving_the_optimisation_problem}
Let us now combine the above approximations and formulate a
tractable optimisation problem, which is amenable to standard
approaches. In addition to searching over the input sequence $\uu$, we
will also treat $\epsilon$ as a slack variable, where the aim is to
reduce $\epsilon$ towards some user defined lower bound $\delta \geq 0$, but
with the possibility that $\epsilon > \delta$ in order that the
constraints form a non-empty feasible set.

Towards this, in what follows it will be convenient to
define a new variable $\zz$ as
\begin{align}
    \zz \triangleq (\uu,\epsilon)
\end{align}
and the $j^{\text{th}}$ element of a new function $g^M \in \R^{n_{cx}}$ as
\begin{align}
g_j^M(\uu,\gamma) &\triangleq \frac{1}{M}\sum_{i=1}^M
       \sigma(c^j_x(x_t^i, u_{t},\theta^i,\uu,\ww^i,),\gamma).
\end{align}
Here, we have dropped the explicit reliance on
$x_t^i, \theta^i$ and $\ww^i$ in order to reduce notation
overhead. The optimisation problem can then be stated as

\begin{small}
  \begin{subequations}
    \label{eq:approxprob}
\begin{align}
    \zz^\star(\gamma) = \arg\,&\min_{\zz}  \quad \eta (\epsilon - \epsilon_0)^2 +  \frac{1}{M}\sum_{i=1}^M V_t(x_t^i, u_{t},\theta^i, \uu,\ww^i),\\
     &\text{s.t.} \hspace{6.mm} c_u(\uu) \succeq 0, \qquad \epsilon \geq \delta, \\ 
     & \hspace{10mm} g^M(\uu,\gamma) \succeq 1-\epsilon.
\end{align}
\end{subequations}
\end{small}

In the above, $\eta > 0$ is a user-defined weighting term and
$\epsilon_0 \leq 0$ is used to offset the penalty on this slack
variable.

The basic approach is to solve a sequence of these problems as
$\gamma \to 0$, thus recovering the indicator function 
constraints in the limit. On the one hand, this is a nonlinear (and
non-convex) programming problem with inequality constraints, and it is
possible to employ standard software for solving this problem. On the
other hand, it is not essential to obtain high precision solutions for
each $\gamma$ value, since small changes in $\gamma$ are likely to
produce small changes in the solution.

Motivated by this, we derive an interior-point method (based on
the log-barrier approach) that is amenable to solving a sequence of
problems as $\gamma \to 0$ (see e.g. \cite{fiacco1990nonlinear,wills2004barrier}). To this end, the log-barrier method
proceeds by first amending the cost with a log-barrier for each
inequality constraint (where the logarithm of a vector is to be interpreted element-wise)

{\small
  \begin{align}
    \notag
    V_t^M(\zz,\gamma,\mu) &\triangleq \eta (\epsilon - \epsilon_0)^2 +  \frac{1}{M}\sum_{i=1}^M V(x_t^i, u_{t},\theta^i, \uu,\ww^i) - \mu \ln(\epsilon - \delta)\\
    \label{eq:log_barrier_cost}
    &-\mu\ln(c_u(\uu)) -\mu \ln \left ( g^M(\uu,\gamma) - 1 + \epsilon \right ).
\end{align}
}%

In the above, $\mu > 0$ is a barrier weighting term, and
importantly, in the limit as $\mu \to 0$, the solution to
{\small
\begin{align}
  \label{eq:unconsol}
  \zz(\gamma,\mu)^\star &= \arg \min_{\zz} V^M_t(\zz,\gamma,\mu) 
\end{align}
}
coincides with the solution to \eqref{eq:approxprob}. The main
benefit of this approach is that problem \eqref{eq:unconsol} is
directly amenable to second-order unconstrained optimisation methods,
such as Newton's method. Note, this requires some care in ensuring
that iterates remain in the feasible region of the log-barrier terms,
which is straighforwardly handled, for example, by a modified line
search procedure.

With this as background, the approach detailed in
Algorithm~\ref{alg:optim_solve} proceeds by solving the sequence of
problems \eqref{eq:unconsol} where the barrier weighting $\mu$ and the
sigmoid function parameter $\gamma$ are gradually reduced to zero.

\begin{algorithm}[!ht]
  \caption{\textsf{Control Action}}
  \footnotesize
  \textsc{Inputs:} $u_t$, $\ww^{1:M}$, $x_t^{1:M}$, $\theta^{1:M}$,
  $\delta$, $c_1$, $\tau_1,\tau_2,\tau_3, \eta_1, \eta_2 > 0$ and max\_iter.\\ \textsc{Output:} $\zz^\star$.
  \algrule[.4pt]
  \begin{algorithmic}[1]
    \STATE Initialise $\mu >> 0$ and $\gamma \approx 1.0$
    \STATE Initialise $\epsilon > \delta$ s.t. $g^M(\uu,\gamma) \geq 1 - \epsilon$
    \STATE Initialise $\uu$ s.t. $c_u(\uu) > 0$
    \FOR{$i < \text{max\_iter}$}
        \STATE Compute the cost $c = V^M(\zz,\gamma,\mu)$ given by \eqref{eq:log_barrier_cost}
        \STATE Compute the gradient $g = \nabla_{\zz}V^M(\zz,\gamma,\mu)$
        \STATE Compute the Hessian $H = \nabla^2_{\zz}V^M(\zz,\gamma,\mu)$
        \STATE Compute the search direction $p = - H^{-1}g$
        \IF{$p^\Transp g > 0$}
        \STATE Find $\lambda > 0$ such that $p = - (H+\lambda
        I)^{-1}g$ ensures that $p^\Transp g < 0$
        \ENDIF
        \STATE Find a step-length $\alpha > 0$ s.t. ${V}^M(\zz + \alpha p, \gamma,\mu) < c + c_1  \alpha  p^{\Transp}g$
        \STATE Set $\zz = \zz +\alpha p$
        \IF{$|p^\Transp g| < \tau_1$ and $\mu \leq \tau_2$ and $\gamma \leq \tau_3$}
        \STATE Terminate with solution $\zz$
        \ELSIF{$|p^\Transp g| < \tau_1$}
          \STATE set $\mu = \max(\eta_1\mu,\tau_2)$
          \STATE set $\gamma = \min(\eta_2\gamma,\tau_3)$
          \STATE determine $\epsilon$ s.t. $g^M(\uu,\gamma) \geq 1 - \epsilon$ 
          \STATE update $\zz$ with the new values of $\epsilon$
        \ENDIF
    \ENDFOR
  \end{algorithmic}
  \label{alg:optim_solve}
\end{algorithm}

\subsection{Resulting Data to Controller} 
\label{sub:on_line_sampling_and_control}
In Algorithm~\ref{alg:smpc} we summarise the resulting solution, which is
demonstrated in Section~\ref{sec:sims}. 
\begin{algorithm}[!ht]
  \caption{\textsf{Data to Controller}}
  \footnotesize
    \textsc{Inputs:} $u_1$
  \algrule[.4pt]
  \begin{algorithmic}[1]
\FOR{$t \geq 1$}
\STATE apply the control action $u_t$;
\STATE measure $y_t$;
\STATE sample $x_t^i, \theta^i,\ww^i \sim p(x_t, \theta, \ww|y_{1:t},
u_{1:t})$ for $i=1,\dots,M$ via HMC;
\STATE compute $\uu^\star$ using Algorithm~\ref{alg:optim_solve};
\STATE set $u_{t+1}$ to the first element of $\uu^*$.
\ENDFOR
\end{algorithmic}
  \label{alg:smpc}
\end{algorithm}






\section{SIMULATIONS}
\label{sec:sims}

\subsection{Pedagogical Example} 
\label{sub:pedagogical_example}
First, the proposed approach is applied to a \red{non-control affine} first-order system to demonstrate the dual accomplishment of learning and control. Additionally, it highlights how the full probabilistic information about the states and parameters informs the control actions. Let us consider
\begin{subequations}
\begin{align}
    x_{t+1} &= a x_t + b \sin(u_t) + w_t, &  w_t&\sim\mathcal{N}(0,q^2), \\
    y_t &= x_t + e_t, & e_t&\sim\mathcal{T}(\nu,0,r^2),
\end{align}
\end{subequations}
\red{where $\mathcal{T}(\nu,0,r^2)$ is a Student's T distribution with $\nu$ degrees of freedom and scale $r$.} For this example, $a=0.9$, $b=0.2$, \red{$r = 0.05$, $\nu=4$}, and $q=0.05$, the constraints $P(0 \leq \xx \leq 1.2) \geq 0.95$ are considered, and a set point of $x^*=1.0$ is chosen.

The system was simulated for $T=50$ discrete time steps. At each time step,  Algorithm~\ref{alg:smpc} was used to sample from the conditional distribution \red{$p(x_t, \theta | y_{1:t}, u_{1:t},\nu)$ where $\theta = \{a, b, q, r\}$ in this case --- with target and achieved acceptance rates of $80\%$ and $79\%$ ---} and to calculate the next control action $u_t$ using a horizon of $N=10$. 
Note that given a sample $q^i$ we can easily sample $\ww^i$.
\red{The computation time per iteration was in the order of $ \SI{4}{\second}$.}

\begin{figure}[!htb]
    \centering
    \includegraphics[width=0.95\linewidth]{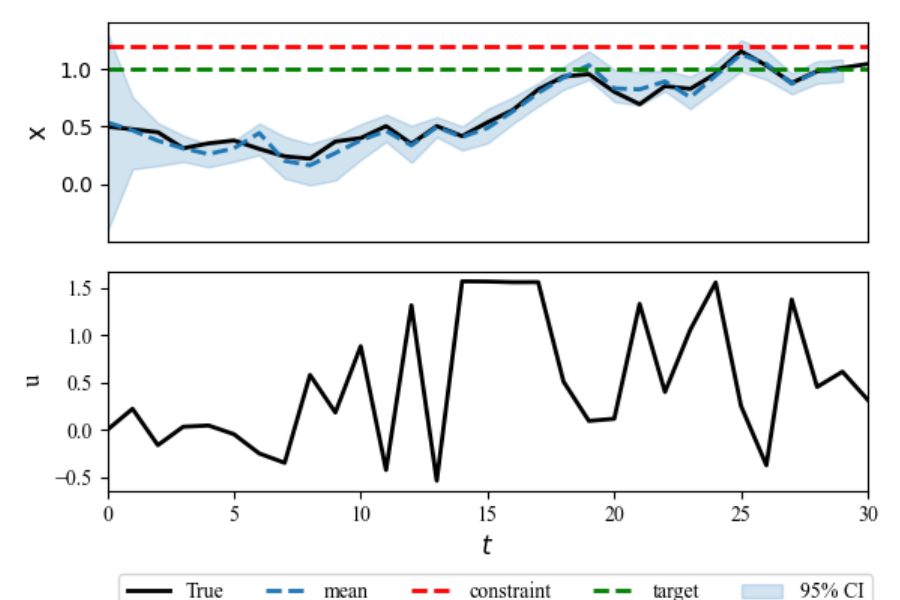}
    \caption{Simulation results of the pedagogical first-order example. Top: the estimated and true state. Bottom: the control action.}
    \label{fig:order1_x_u}
\end{figure}

\begin{figure}[!htb]
    \centering
    \includegraphics[width=0.95\linewidth]{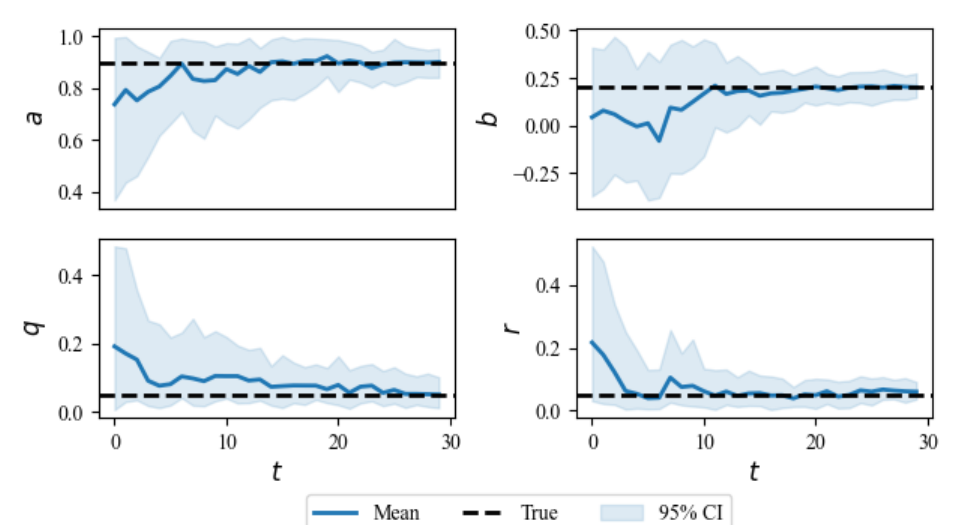}
    \caption{Evolution of the parameter estimates over the simulation.
    }
    \label{fig:order1_params}
\end{figure}

Figure~\ref{fig:order1_x_u} shows the simulated and estimated state along with the control input and Figure~\ref{fig:order1_params} shows the parameter estimates at each time step.
Relatively uninformative priors were placed on the initial state and on the parameters. As such, we can see that the estimates of the parameters and state become more certain as the amount of data $\{y_{1:t},u_{1:t}\}$ increases. Consequently, the control action becomes more aggressive and the state is driven closer whilst still satisfying the chance state constraint. 

\begin{figure}[!htb]
    \centering
    \includegraphics[width=0.95\linewidth]{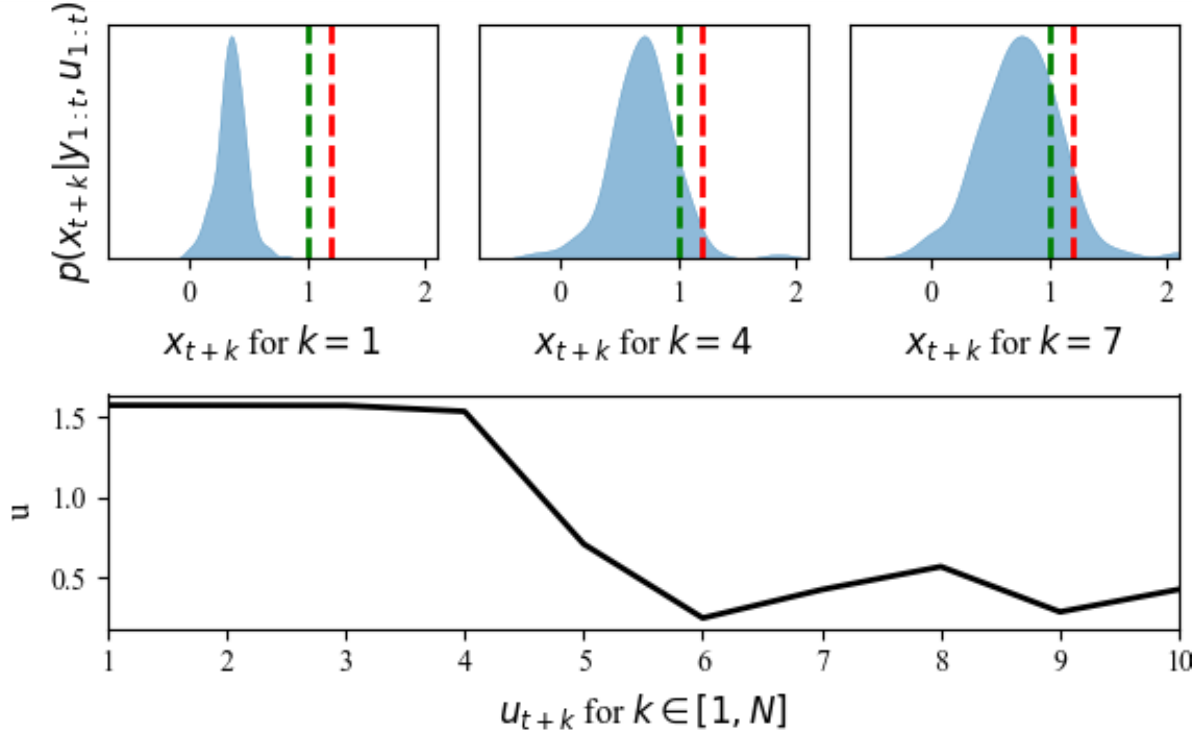}
    \caption{Predicted states $\xx$ and control action $\uu$ over the horizon for a single time step at $t=6$. Top: Example distributions of the predicted state at several discrete time steps over the horizon. Also shown is the set point (dashed green), and state constraint (dashed red). 
    Bottom: the control action over the prediction horizon.}
    \label{fig:order1_mpc_horizon}
\end{figure}

Additionally, examples of the predicted future states, $\xx$, and optimised control actions over the horizon, $\uu$, are provided in Figure~\ref{fig:order1_mpc_horizon}. A trade-off exists between achieving the set point and satisfying the constraint with the desired probability. As a consequence, we can observe that the uncertainty on $\xx$ grows the further into the future the prediction is and as a consequence the forecast control actions become more conservative in order to ensure the desired $95\%$ state constraint satisfaction at these future time steps. 

\subsection{Rotary Inverted Pendulum} 
\label{sub:rotary_inverted_pendulum}
The proposed method is now applied to a simulated rotational inverted pendulum, or Furata pendulum \cite{furuta1992swing}. For the purpose of simulation the parameter values specified for the QUANSER\textsuperscript{\textcopyright}  QUBE-Servo 2 Rotary inverted pendulum in \cite{guides2003laboratories} are used.

Let the state vector used to model the Furata pendulum be $x = \pmat{\vartheta & \alpha & \dot{\vartheta} & \dot{\alpha}}^{\Transp}$, where $\vartheta$ and $\alpha$ are the base arm and pendulum angles respectively, and the controllable input to the system be the motor voltage $V_{m}$.
Then, the continuous time dynamics are given by

\small
\begin{equation}\label{eq:pend_process}
\begin{split}
    &M(\alpha)\pmat{\ddot{\vartheta} \\ \ddot{\alpha}} + \nu(\dot{\vartheta},\dot{\alpha})\pmat{\dot{\vartheta} \\ \dot{\alpha}} = \pmat{\frac{k_m(V_m - k_m\dot{\vartheta})}{R_m} - D_r(\dot{\vartheta}) \\ -\frac{1}{2}m_pL_pg\sin(\alpha)-D_p\dot{\alpha}},\\
         &M(\alpha) \hspace{-0.75mm}=\hspace{-0.75mm} \pmat{m_pL_r^2\hspace{-0.5mm}+\hspace{-0.5mm}\frac{1}{4}m_pL_p^2(1\hspace{-0.5mm}-\hspace{-0.5mm}\cos(\alpha)^2)\hspace{-0.5mm}+\hspace{-0.5mm}J_r & \hspace{-2mm} \frac{1}{2}m_pL_pL_r\cos(\alpha) \\ \frac{1}{2}m_pL_pL_r\cos(\alpha) & J_p + \frac{1}{4}m_pL_p^2}, \\
    &\nu(\dot{\vartheta},\dot{\alpha}) \hspace{-0.75mm}=\hspace{-0.75mm} \pmat{\frac{1}{2}m_pL_p^2\sin(\alpha)\cos(\alpha)\dot{\alpha} & \hspace{-2mm} - \frac{1}{2}m_pL_pL_r\sin(\alpha)\dot{\alpha} \\
    - \frac{1}{4}m_pL_p^2\cos(\alpha)\sin(\alpha)\dot{\vartheta} & 0}, \\
\end{split}
\end{equation}
\normalsize
where $m_p$ is the pendulum mass, $L_r$, $L_p$ and the rod and pendulum lengths, $J_r$, $J_p$ are the rod and pendulum inertias, $R_m$ and $k_m$ are the motor resistance and constant, $D_p$ is the pendulum damping constant, and $D_r(\dot\vartheta)$ is the arm damping function. \red{Since a perfect damping model is unlikely, $D_r(\dot\vartheta) = D_{r0}\text{sign}(\dot\vartheta) + D_{r1}\dot\vartheta + D_{r2}\dot\vartheta^2$ was used for simulation and $D_r(\dot\vartheta)=D_{r1}\dot\vartheta$ for estimation and control.}


The process model used for simulation, estimation and control consists of a 4\textsuperscript{th} order Runge-Kutta discretisation of these continuous time dynamics over a $25$ms sampling time, subsequently disturbed by noise $w_t$. The input voltage from $t$ to $t+1$ is given by a zero-order hold of the control $u_t$.

Measurements are from encoders on the arm and pendulum angle and current measurements from the motor. The resulting measurement model is $y_t = \pmat{\vartheta & \alpha & \frac{V_m-k_m\dot{\vartheta}}{R_m}} + e_t$.
It is assumed that the elements of the process and measurement noise are all independent: \red{$w_t\sim\mathcal{N}(0,\Sigma_w)$ and $e_t\sim\mathcal{N}(0,\Sigma_e)$, i.e. 
 $\Sigma_w = \text{diag}([\num{3e-4}, \num{1e-4}, \num{0.013}, \num{0.013}])^2$ and $\Sigma_e = \text{diag}([\num{1.1e-3}, \num{1e-3}, \num{0.175}])^2$. 
  $u_{t}$ is constrained to  $[-18,18]$, $\vartheta$ is constrained to $\pm 270^\circ$ with $95\%$ probability, and the state cost $\vartheta^2 + 900(\cos(\alpha) - 1)^2+\num{1e-10}\dot\vartheta^2+\num{1e-10}\dot\alpha$ is used.}

The system was simulated for $T=50$ discrete time steps.
At each time step Algorithm~\ref{alg:smpc} samples from the joint conditional distribution of the state $x_t$
and the parameters, $\theta=\{J_r, J_p, K_m, R_m, D_p, D_r, \Sigma_w,\Sigma_e\}$, given by $p(x_t,\theta|y_{1:t})$ \red{--- with target and achieved acceptance rates of $85\%$ and $83\%$ ---} and calculate the control $u_{t+1}$ with a horizon of $N=25$. \red{The computation time per iteration was in the order of $\SI{200}{\second}$.}

\begin{figure}[tb]
    \centering
    \includegraphics[width=0.95\linewidth]{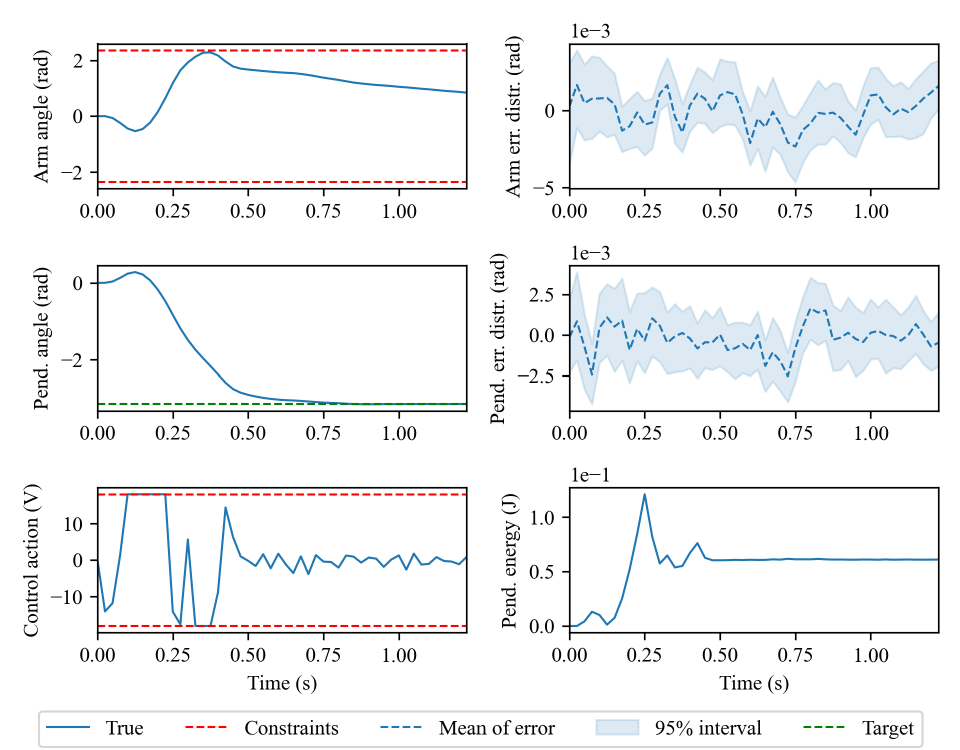}
    \caption{Simulation results for the rotary inverted pendulum. Left: the true states. Right: the distribution of the error between the samples and the true states. Bottom: control action and pendulum arm total energy.}
    \label{fig:inv_pend_state_traj}
\end{figure}

\begin{figure}[tb]
    \centering
    \includegraphics[width=0.95\linewidth]{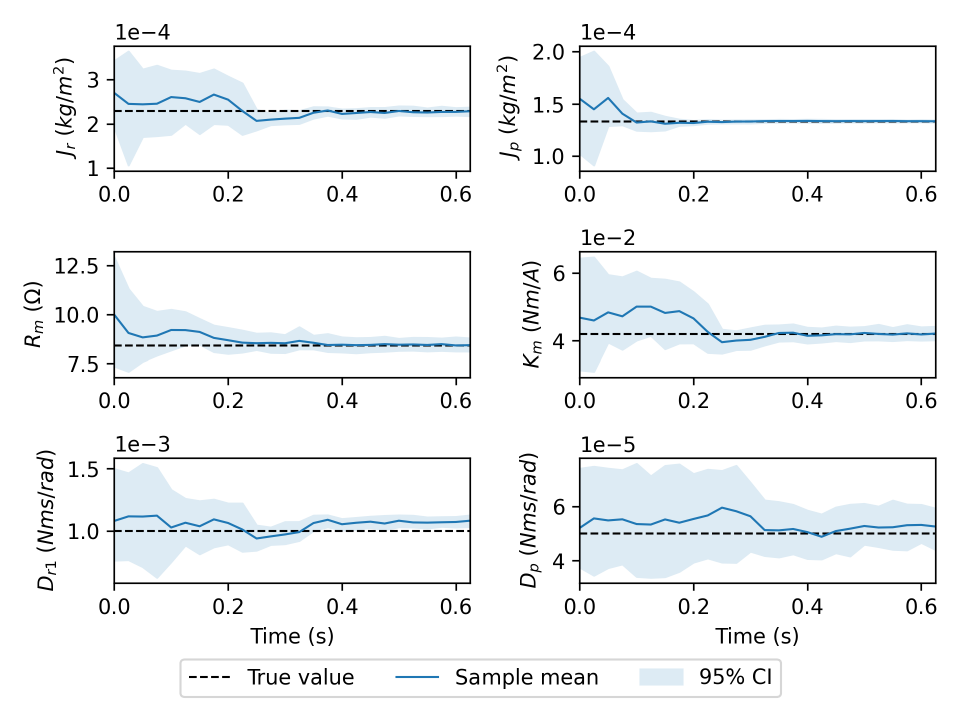}
    \caption{Evolution of the pendulum model parameter estimates over the simulation. Only half the simulation time is shown as the estimates have mostly converged by this stage.}
    \label{fig:pend_params_est}
\end{figure}

Figure~\ref{fig:inv_pend_state_traj} shows the true and estimated arm and pendulum angles as well as the control action and the total energy of the pendulum. Figure~\ref{fig:pend_params_est} shows the model parameter estimates at each discrete time step. These results show that the proposed method is capable of updating the state and parameter estimates as more data becomes available and using the full probabilistic information to determine the control in a nonlinear setting.
An animation illustrating these results is available at \url{https://youtu.be/E3_rUWejXEc}. 

\section{Conclusion and Future Work} 
\label{sec:conclusion_and_future_works}
This paper presents an approach to determine an optimal control action based on observed data.
The problem formulation directly relates existing system data and its inherent uncertainty to future control actions while handling input and state-based constraints: using a chance constraint formulation for the latter. 
The problem is approximated using HMC in order to produce a tractable optimisation problem. We propose a barrier function based approach for solving this problem and demonstrate the utility of the approach on two simulations examples.

It is important to highlight that the proposed approach will ultimately become intractable as the data length grows. Strategies for using only a finite window of past data and updating the prior distributions are therefore an important research avenue to explore. \red{At present, the approach has a significant computational cost which would restrict its practical applications to systems with slower timescales. Future research will focus on reducing this burden by modifying alternate optimisation procedures to address the case of a sequence as $\gamma\to0$ and on developing more more efficient HMC algorithms for state space models.}
Finally, we have not explicitly treated the problem of exploration. This is a key subject within the nonlinear dual control literature, and this deserves further attention.


\bibliographystyle{plain}
\bibliography{references.bib}
\end{document}